\documentclass[11pt]{article}%
\usepackage[dvips]{graphicx}
\usepackage{theorem}
\usepackage{amsfonts}
\usepackage{amssymb}
\usepackage{amsmath}
\usepackage{amsbsy,amsfonts}
\usepackage{amscd}
\usepackage{calc}
\usepackage{color}
\usepackage{longtable}%
\DeclareTextSymbol{\degre}{OT1}{23}
\newcommand{\nc}{\newcommand}
\nc{\rnc}{\renewcommand}

\newcommand{\qed}{\nobreak \ifvmode \relax \else
\ifdim\lastskip<1.5em \hskip-\lastskip
\hskip1.5em plus0em minus0.5em \fi \nobreak
\vrule height0.75em width0.5em depth0.25em\fi}
\nc{\dsp}{\displaystyle}
\rnc{\vector}{\boldsymbol} \nc{\Grad}{\vector{\mathrm{grad}} \; }
\nc{\Span}{\text{span}} \nc{\Rr}{\mathbb{R}}
\nc{\Rtwo}{\Rr^{2}} \nc{\Rthree}{\Rr^{3}}
\nc{\Pp}{\mathbb{P}}
\rnc{\thefootnote}{\fnsymbol{footnote}}
\begin{document}

{\LARGE
}

\begin{center}
{\LARGE
\textbf{Divergence-free $\mathbf{\mathcal{H}}(\mathbf{div})$-conforming
hierarchical bases for magnetohydrodynamics (MHD)}}
\end{center}

{\LARGE
}


\begin{center}
Wei Cai$^{\dagger}$\footnotemark[1], Jian Wu and Jianguo Xin

{\footnotesize {\emph{$^{\dagger}$Department of Mathematics and Statistics,
University of North Carolina at Charlotte, Charlotte, NC 28223, USA.}} }
\vspace*{.2cm}

{\footnotesize \footnotetext[1]{{\footnotesize Corresponding author. Tel.:
+1-704-687-4581; fax: +1-704-687-6415. E-mail address: wcai@uncc.edu (W.
Cai).}}
}
\end{center}

\begin{abstract}
In order to solve the magnetohydrodynamics (MHD) equations with a
$\mathbf{\mathcal{H}}(\mathbf{div})$-conforming element, a novel approach is
proposed to ensure the exact divergence-free condition on the magnetic field.
The idea is to add on each element an extra interior bubble function from
higher order hierarchical $\mathbf{\mathcal{H}}(\mathbf{div})$-conforming
basis. Four such hierarchical bases for the $\mathbf{\mathcal{H}}%
(\mathbf{div})$-conforming quadrilateral, triangular, hexahedral and
tetrahedral elements are either proposed (in the case of tetrahedral) or
reviewed. Numerical results have been presented to show the linear
independence of the basis functions for the two simplicial elements. Good
matrix conditioning has been confirmed numerically up to the fourth order for
the triangular element and up to the third order for the tetrahedral element.

\vspace*{.2cm}

{\noindent\textit{Keywords}: Hierarchical bases, $\mathbf{\mathcal{H}%
}(\mathbf{div})$-conforming elements, Divergence-free condition}

\end{abstract}

\rnc{\thefootnote}{\arabic{footnote}} \setcounter{footnote}{0}



\section{Introduction}

The magnetohydrodynamics (MHD) equations describe the dynamics of a charged
system under the interaction with a magnetic field and the conservation of the
mass, momentum and energy for the plasma system. Such a dynamics is considered
constrained as the magnetic field of the system is evolved with the constraint
of zero divergence, namely, $\nabla\cdot\mathbf{B}=0$. Numerical modeling of
plasmas has shown that the observance of the zero divergence of the magnetic
field plays an important role in reproducing the correct physics in the plasma
fluid \cite{Brackbill80}. Various numerical techniques have been devised to
ensure the computed magnetic field to maintain divergence-free \cite{Toth00}.
In the original work of \cite{Brackbill80} a projection approach was used to
correct the magnetic field to have a zero divergence.

A more natural way to satisfy this constraint is through a class of the
so-called constrained transport (CT) numerical methods based on the ideas in
\cite{Evans88}. As noted in \cite{Monk03}, a piecewise $\mathbf{\mathcal{H}%
}(\mathbf{div})$ vector field on a finite element triangulation of a spatial
domain can be a global $\mathbf{\mathcal{H}}(\mathbf{div})$ field if and only
if the normal components on the interface of adjacent elements are continuous.
Thus, in most of the CT algorithms for the MHD, the surface averaged magnetic
flux over the surfaces of a 3-D element is used to represent the magnetic
field while the volume averaged conserved quantities (mass, momentum, and
energy) are used.

In the two seminal papers \cite{Ned1,Ned2}, N\'{e}d\'{e}lec proposed to use
quantities (moments of tangential components of vector fields) on edges and
faces to define the finite dimensional space in $\mathbf{\mathcal{H}%
}(\mathbf{div})$ and $\mathbf{\mathcal{H}}(\mathbf{curl})$. The specific
construction of the basis functions in both spaces, specifically in
$\mathbf{\mathcal{H}}(\mathbf{div})$, can be done in various ways such as the
hierarchical type of basis proposed in \cite{AC03B} for both
$\mathbf{\mathcal{H}}(\mathbf{div})$ and $\mathbf{\mathcal{H}}(\mathbf{curl}%
)$. Unfortunately, the proposed hierarchical basis for $\mathbf{\mathcal{H}%
}(\mathbf{div})$ in \cite{AC03B} turns out to be erroneous as it can be easily
checked that for quadratic polynomial approximation the proposed edge-based
basis functions happens to be linearly dependent.

In this paper, we will first present a new hierarchical basis functions in
$\mathbf{\mathcal{H}}(\mathbf{div})$ in tetrahedral, for completeness,
together with a review of the hierarchical basis of $\mathbf{\mathcal{H}%
}(\mathbf{div})$ for rectangles in 2-D and for cube in 3-D. An important
common feature of the hierarchical basis functions is the fact that for order
$p\geq3$ in the case of simplexes, the basis functions will include interior
bubble basis functions which have zero normal components on the whole boundary
of the element boundaries. Therefore, a simple way to enforce the
divergence-free on each element can be easily accomplished by adding one
single $(p+1)$-th order (any $q$-th order, $q\geq\max(m,p+1),m=4$ for
tetrahedral element, $m=3$ for triangular element, and $m=2$ for quadrilateral
or hexahedral element) interior bubble function to a $p$-th order
$\mathbf{\mathcal{H}}(\mathbf{div})$ basis. This extra bubble basis will be
able to satisfy the local divergence-free condition. Due to the fact of normal
continuity of the $p$-th order basis across the element interface and zero
normal component of the added-in $p+1$-th order interior bubble functions, the
augmented function space satisfies divergence-free condition globally.
$\mathbf{\mathcal{H}}(\mathbf{div})$ basis on general meshes other than the
reference elements mentioned above are usually constructed by a Piola
transform  \cite{falk11} and more recent studies of the $\mathbf{\mathcal{H}}(\mathbf{div})$  basis
functions on general quadrilateral and hexahedral elements can be found in
\cite{falk11} \cite{kwak11} \cite{wheeler12}, however, they will not be
discussed further in this paper.

The rest of the paper is organized as follows. The constructions of the
locally divergence-free bases are given in Section 2-5 for rectangular and
triangular elements in 2-D, and cubic and tetrahedral elements in 3-D. The
divergence-free condition is discussed in Section 6. Numerical results on the
matrix conditioning are given in Section 7. Concluding remarks are given in
Section 8.

\section{Basis functions for the quadrilateral element}

\label{sec:construct2dquad} In \cite{Zag} Zaglmayr gave a hierarchical basis
for quadrilateral $\mathbf{\mathcal{H}}(\mathbf{div})$-conforming element. In
this section we summarize the result in \cite{Zag}.

The basis functions are constructed on the reference element - a unit square
$\mathcal{Q}:=[0,1]^{2}$ \cite{Zag} with vertexes of $V_{1}(0,0)$,
$V_{2}(1,0)$, $V_{3}(1,1)$ and $V_{4}(0,1)$. The coordinate system for the
reference element is in terms of the variables $(\xi,\eta)$. A bilinear
function $\lambda_{i}$ which is associated with a specific vertex $V_{i}$ has
been utilized for the construction. The multiplicative factors of the bilinear
function $\lambda_{i}$ have been used to form the linear function $\sigma_{i}%
$, viz.
%

\[
\lambda_{1}:=(1-\xi)(1-\eta),\quad\sigma_{1}:=(1-\xi)+(1-\eta),
\]%
\[
\lambda_{2}:=\xi(1-\eta),\quad\sigma_{2}:=\xi+(1-\eta),
\]%
\[
\lambda_{3}:=\xi\eta,\quad\sigma_{3}:=\xi+\eta,
\]%
\begin{equation}
\lambda_{4}:=(1-\xi)\eta,\quad\sigma_{4}:=(1-\xi)+\eta. \label{eqn:bi4}%
\end{equation}
The bilinear function has this favorable property%

\begin{equation}
\label{eqn:bilin}\lambda_{i}|_{V_{j}} = \delta_{ij},
\end{equation}
where $\delta_{ij}$ is the Kronecker delta. The edge $e := [V_{i}, V_{j}]$
which points from vertex $V_{i}$ to vertex $V_{j}$ is parameterized by
\begin{equation}
\label{eqn:edgepar}\zeta_{e} := \sigma_{j} - \sigma_{i} \in[-1, 1].
\end{equation}
For convenience of basis construction, the linear edge-extension parameter is
also defined, viz.
\begin{equation}
\label{eqn:edgext}\lambda_{e} := \lambda_{i} + \lambda_{j} \in[0, 1],
\end{equation}
which is one on edge $e := [V_{i}, V_{j}]$ and zero on the opposite edge. Note
that the unit tangential vector $\tau_{e}$ and the outward unit normal vector
$\mathbf{n}_{e}$ can be deduced as
\begin{equation}
\label{eqn:edgededu}\vec{\tau}_{e} = \frac{1}{2} \nabla\zeta_{e},
\quad\mathbf{n}_{e} = \nabla\lambda_{e}.
\end{equation}

\subsection{Edge-based functions}

\label{sec:face}

These functions are further grouped into two categories: the lowest-order and
higher-order functions.

\bigskip

\noindent{\underline{Lowest-order functions}:}

\bigskip

These functions are associated with the four edges. By construction each
function is perpendicular to the associated edge and has unit normal component
on the associated edge. Furthermore the divergence of each function is unit.
The shape function is given by
\begin{equation}
\label{eqn:edgeface11}\psi_{e_{i}}^{RT_{0}} = \frac{1}{2} \lambda_{e_{i}}
\left(  \nabla\times\zeta_{e_{i}}\right)  , \quad i = \{1, 2, 3, 4\},
\end{equation}
which has the property
\begin{equation}
\label{eqn:edgelow1}\vec{\tau}_{e_{i}} \cdot\psi_{e_{i}}^{RT_{0}} = 0,
\quad\mathbf{n}_{e_{i}} \cdot\psi_{e_{i}}^{RT_{0}}|_{e_{i}} = 1, \quad
\nabla\cdot\psi_{e_{i}}^{RT_{0}} = 1.
\end{equation}

\bigskip

\noindent{\underline{Higher-order functions}:}

\bigskip

The function is taken to be a curl field of a scalar function in order to be
free of divergence. The basis function is
\begin{equation}
\label{eqn:edgehigher1}\psi_{e_{i}}^{j+1} = \nabla\times\left(  \lambda
_{e_{i}} L_{j+2} \left(  \zeta_{e_{i}}\right)  \right)  , \quad i = \{1, 2, 3,
4\}, \quad0 \le j \le p - 1,
\end{equation}
where the function $L_{n}(\bullet)$ is the so-called integrated Legendre
polynomial of degree $n$ \cite{Zag}. The obvious property is
\begin{equation}
\label{eqn:facebubbleorthonormal1}\nabla\cdot\psi_{e_{i}}^{j+1} = 0, \quad i =
\{1, 2, 3, 4\}, \quad0 \le j \le p - 1.
\end{equation}

\subsection{Interior functions}

\label{sec:interior} The interior functions are classified into three categories.

\bigskip

\noindent{\underline{Type 1}: (curl field)}

\bigskip

The shape functions are given as
\begin{equation}
\label{eqn:type1}\psi_{ij}^{Q_{1}} = \nabla\times\left(  L_{i+2}(2\xi-
1)L_{j+2}(2\eta- 1)\right)  , \quad0 \le i, j \le p - 1.
\end{equation}
These functions are divergence free, viz.
\begin{equation}
\label{eqn:type1pro}\nabla\cdot\psi_{ij}^{Q_{1}} = 0, \quad0 \le i, j \le p -
1.
\end{equation}

\bigskip

\noindent{\underline{Type 2}:}

\bigskip

The formula of these functions is given as
\begin{equation}
\label{eqn:type2}\psi_{ij}^{Q_{2}} = L_{i+2}(2\xi- 1)L_{j+2}^{\prime}(2\eta-
1) \hat{\mathbf{i}}_{\xi} + L_{i+2}^{\prime}(2\xi- 1)L_{j+2}(2\eta- 1)
\hat{\mathbf{j}}_{\eta}, \quad0 \le i, j \le p - 1.
\end{equation}
On the boundary of the reference element the normal component of these
functions vanishes, viz.
\begin{equation}
\label{eqn:type2pro}\mathbf{n}_{e} \cdot\psi_{ij}^{Q_{2}}|_{e} = 0, \quad0 \le
i, j \le p - 1,
\end{equation}
which justifies to be one type of interior functions.

\bigskip

\noindent{\underline{Type 3}:}

\bigskip

The formula of these functions is
\begin{equation}
\psi_{i}^{Q_{3}^{\xi}}=L_{i+2}(2\xi-1)\hat{\mathbf{i}}_{\xi},\quad\psi
_{i}^{Q_{3}^{\eta}}=L_{i+2}(2\eta-1)\hat{\mathbf{j}}_{\eta},\quad0\leq i\leq
p-1. \label{eqn:type3}%
\end{equation}
Again on the boundary of the reference element the normal component of these
functions vanishes, viz.
\begin{equation}
\mathbf{n}_{e}\cdot\psi_{i}^{Q_{3}^{\xi}}|_{e}=0,\quad\mathbf{n}_{e}\cdot
\psi_{i}^{Q_{3}^{\eta}}|_{e}=0,\quad0\leq i\leq p-1. \label{eqn:type3pro}%
\end{equation}

Table 1 shows the decomposition of the space $Q_{p+1,p}\times Q_{p,p+1}$ for
the $\mathbf{\mathcal{H}}(\mathbf{div})$-conforming quadrilateral
element.\begin{table}[th]
\begin{center}%
\begin{tabular}
[c]{|c|c|}\hline
Decomposition & Dimension\\\hline
Edge-based functions & $4(p+1)$\\\hline
Interior functions & $2p(p+1)$\\\hline
Total & $2(p+2)(p+1)=Dim(Q_{p+1,p}\times Q_{p,p+1})$\\\hline
\end{tabular}
\end{center}
\caption{Decomposition of the space $Q_{p+1,p}\times Q_{p,p+1}$ for the
$\mathbf{\mathcal{H}}(\mathbf{div})$-conforming quadrilateral element.}%
\end{table}

\section{Basis functions for the hexahedral element}

\label{sec:construct3hex} The reference element is defined for a unit cube
$\mathcal{H}:=[0,1]^{3}$ in \cite{Zag}. The vertexes of the cube are
$V_{1}(0,0,0)$, $V_{2}(1,0,0)$, $V_{3}(1,1,0)$, $V_{4}(0,1,0)$, $V_{5}%
(0,0,1)$, $V_{6}(1,0,1)$, $V_{7}(1,1,1)$, and $V_{8}(0,1,1)$. The basis
functions are expressed in terms of the trilinear function $\lambda_{j}$,
which is one at the vertex $V_{j}$ and zero at all other vertexes. Along with
the linear function $\sigma_{j}$ and in terms of the coordinates variables
$(\xi,\eta,\zeta)$ they are given as
%

\[
\lambda_{1}:=(1-\xi)(1-\eta)(1-\zeta),\quad\sigma_{1}:=(1-\xi)+(1-\eta
)+(1-\zeta),
\]%
\[
\lambda_{2}:=\xi(1-\eta)(1-\zeta),\quad\sigma_{2}:=\xi+(1-\eta)+(1-\zeta),
\]%
\[
\lambda_{3}:=\xi\eta(1-\zeta),\quad\sigma_{3}:=\xi+\eta+(1-\zeta),
\]%
\[
\lambda_{4}:=(1-\xi)\eta(1-\zeta),\quad\sigma_{4}:=(1-\xi)+\eta+(1-\zeta).
\]%
\[
\lambda_{5}:=(1-\xi)(1-\eta)\zeta,\quad\sigma_{5}:=(1-\xi)+(1-\eta)+\zeta,
\]%
\[
\lambda_{6}:=\xi(1-\eta)\zeta,\quad\sigma_{6}:=\xi+(1-\eta)+\zeta,
\]%
\[
\lambda_{7}:=\xi\eta\zeta,\quad\sigma_{7}:=\xi+\eta+\zeta,
\]%
\begin{equation}
\lambda_{8}:=(1-\xi)\eta\zeta,\quad\sigma_{8}:=(1-\xi)+\eta+\zeta.
\label{eqn:bi38}%
\end{equation}
The edge $e:=[V_{i},V_{j}]$ which points from vertex $V_{i}$ to vertex $V_{j}$
is parameterized by%

\begin{equation}
\label{eqn:edgepar3}\mu_{e} := \sigma_{j} - \sigma_{i} \in[-1, 1].
\end{equation}
The tangential vector associated with the edge $e$ is given by $\vec{\tau}_{e}
= \frac{1}{2}\nabla\mu_{e}$. The edge extension parameter $\lambda_{e} :=
\lambda_{i} + \lambda_{j} \in[0, 1]$ is one on the edge $e$ and zero on all
other edges that are parallel to the edge $e$. The face $f = [V_{i}, V_{j},
V_{k}, V_{\ell}]$ where the vertexes $V_{i}$ and $V_{k}$ are not connected by
an edge can be parameterized by
\begin{equation}
\label{eqn:edgepara3}(\xi_{f}, \eta_{f}) := (\sigma_{i} - \sigma_{j},
\sigma_{i} - \sigma_{\ell}) \in[-1,1] \times[-1,1].
\end{equation}
The linear face extension parameter $\lambda_{f} = \lambda_{i} + \lambda_{j} +
\lambda_{k} + \lambda_{\ell}$ is equal to one on the face $f$ and zero on the
opposite face. The outward unit normal vector of face $f$ can be obtained by
$\mathbf{n}_{f} = \nabla\lambda_{f}$.

\subsection{Face-based functions}

\label{sec:face3Dhex} In this subsection we record the results in \cite{Zag}.
We have also fixed one error in \cite{Zag}. These functions are associated
with the six faces whose formulas are classified into two groups.

\bigskip

\noindent{\underline{Lowest-order Raviart-Thomas functions}:}%

\begin{equation}
\label{eqn:RT0}\psi_{f_{i}}^{\mathcal{RT}_{0}} = \lambda_{f} \mathbf{n}_{f},
\quad i = 1, 2, \cdots, 6.
\end{equation}

\bigskip

\noindent{\underline{Higher-order functions}: (divergence-free)}

\bigskip

These functions are constructed as curl fields of certain vectors in order to
be divergence-free. The formulas are given as \cite{Zag}
%

\begin{equation}
\psi_{i,j}^{f_{k}} = \nabla\times\left(  \lambda_{f} \left(  L_{j+2}(\eta_{f})
\nabla L_{i+2}(\xi_{f}) - L_{i+2}(\xi_{f}) \nabla L_{j+2}(\eta_{f})\right)
\right)  , \;\; 0 \le i, j \le p-1, \;\; k = 1,2,\cdots,6.
\end{equation}
\begin{equation}
\label{eqn:facehigherhexi}\psi_{i}^{f_{k}} = \nabla\times\left(  \lambda_{f}
L_{i+2}(\xi_{f})\nabla\eta_{f}\right)  , \quad0 \le i \le p-1, \quad k =
1,2,\cdots,6.
\end{equation}
\begin{equation}
\label{eqn:facehigherhexj}\psi_{j}^{f_{k}} = \nabla\times\left(  \lambda_{f}
L_{j+2}(\eta_{f})\nabla\xi_{f}\right)  , \quad0 \le j \le p-1, \quad k =
1,2,\cdots,6.
\end{equation}

\subsection{Interior functions}

The interior functions are further classified into three categories. The
triplet $(\xi_{1}, \eta_{2}, \zeta_{3}) := (2\xi-1, 2\eta- 1, 2\zeta- 1)$ is
used in the formulas. The function $\ell_{n}(\bullet)$ is the classical
un-normalized Legendre polynomial of degree $n$. While we here record the
results in \cite{Zag}, we have implemented the correction of a number of
mistakes in \cite{Zag} as well.

\bigskip

\noindent{\underline{Type 1}: (divergence-free)}

\bigskip

These functions are taken to be the curl fields of certain vector functions
that are associated with the $\mathbf{\mathcal{H}}(\mathbf{curl})$-conforming element.

%

\[
\psi_{i,j,k}^{\mathcal{H}_{1}^{1}}=4L_{i+2}(\xi_{1})\ell_{j+1}(\eta_{2}%
)\ell_{k+1}(\zeta_{3})\hat{\mathbf{i}}_{\xi}-4\ell_{i+1}(\xi_{1})\ell
_{j+1}(\eta_{2})L_{k+2}(\zeta_{3})\hat{\mathbf{k}}_{\zeta},\quad0\leq
i,j,k\leq p-1.
\]%
\[
\psi_{i,j,k}^{\mathcal{H}_{1}^{2}}=4\ell_{i+1}(\xi_{1})L_{j+2}(\eta_{2}%
)\ell_{k+1}(\zeta_{3})\hat{\mathbf{j}}_{\eta}-4\ell_{i+1}(\xi_{1})\ell
_{j+1}(\eta_{2})L_{k+2}(\zeta_{3})\hat{\mathbf{k}}_{\zeta},\quad0\leq
i,j,k\leq p-1.
\]%
\[
\psi_{j,k}^{\mathcal{H}_{1}^{3}}=2\ell_{j+1}(\eta_{2})L_{k+2}(\zeta_{3}%
)\hat{\mathbf{k}}_{\zeta}-2L_{j+2}(\eta_{2})\ell_{k+1}(\zeta_{3}%
)\hat{\mathbf{j}}_{\eta},\quad0\leq j,k\leq p-1.
\]%
\[
\psi_{i,k}^{\mathcal{H}_{1}^{4}}=2L_{i+2}(\xi_{1})\ell_{k+1}(\zeta_{3}%
)\hat{\mathbf{i}}_{\xi}-2\ell_{i+1}(\xi_{1})L_{k+2}(\zeta_{3})\hat{\mathbf{k}%
}_{\zeta},\quad0\leq i,k\leq p-1.
\]%
\begin{equation}
\psi_{i,j}^{\mathcal{H}_{1}^{5}}=2L_{i+2}(\xi_{1})\ell_{j+1}(\eta_{2}%
)\hat{\mathbf{i}}_{\xi}-2\ell_{i+1}(\xi_{1})L_{j+2}(\eta_{2})\hat{\mathbf{j}%
}_{\eta},\quad0\leq i,j\leq p-1. \label{eqn:inter3dT15}%
\end{equation}

\noindent{\underline{Type 2}:}

\bigskip

These functions are linear combinations of certain components in the above type.

%

\[
\psi_{i,j,k}^{\mathcal{H}_{2}^{1}}=L_{i+2}(\xi_{1})\ell_{j+1}(\eta_{2}%
)\ell_{k+1}(\zeta_{3})\hat{\mathbf{i}}_{\xi}+\ell_{i+1}(\xi_{1})L_{j+2}%
(\eta_{2})\ell_{k+1}(\zeta_{3})\hat{\mathbf{j}}_{\eta},\quad0\leq i,j,k\leq
p-1.
\]%
\[
\psi_{j,k}^{\mathcal{H}_{2}^{2}}=L_{j+2}(\eta_{2})\ell_{k+1}(\zeta_{3}%
)\hat{\mathbf{j}}_{\eta}+\ell_{j+1}(\eta_{2})L_{k+2}(\zeta_{3})\hat
{\mathbf{k}}_{\zeta},\quad0\leq j,k\leq p-1.
\]%
\[
\psi_{i,k}^{\mathcal{H}_{2}^{3}}=\ell_{i+1}(\xi_{1})L_{k+2}(\zeta_{3}%
)\hat{\mathbf{k}}_{\zeta}+L_{i+2}(\xi_{1})\ell_{k+1}(\zeta_{3})\hat
{\mathbf{i}}_{\xi},\quad0\leq i,k\leq p-1.
\]%
\begin{equation}
\psi_{i,j}^{\mathcal{H}_{2}^{4}}=L_{i+2}(\xi_{1})\ell_{j+1}(\eta_{2}%
)\hat{\mathbf{i}}_{\xi}+\ell_{i+1}(\xi_{1})L_{j+2}(\eta_{2})\hat{\mathbf{j}%
}_{\eta},\quad0\leq i,j\leq p-1. \label{eqn:inter3dT24}%
\end{equation}

\noindent{\underline{Type 3}:}

\bigskip

These functions are taken as certain components in Type 2.

%

\[
\psi_{i}^{\mathcal{H}_{3}^{1}}=L_{i+2}(\xi_{1})\hat{\mathbf{i}}_{\xi}%
,\quad0\leq i\leq p-1.
\]%
\[
\psi_{j}^{\mathcal{H}_{3}^{2}}=L_{j+2}(\eta_{2})\hat{\mathbf{j}}_{\eta}%
,\quad0\leq j\leq p-1.
\]%
\begin{equation}
\psi_{k}^{\mathcal{H}_{3}^{3}}=L_{k+2}(\zeta_{3})\hat{\mathbf{k}}_{\zeta
},\quad0\leq k\leq p-1. \label{eqn:inter3dT33}%
\end{equation}

Two remarks are in place.

\begin{itemize}
\item All the interior basis functions are linearly independent, which can be
verified easily.

\item The normal traces of these interior functions vanish on the boundary
$\partial\mathcal{H}$ of the reference hexahedral element. This is due to the
fact that on a certain face $f$ either one of the standard unit vectors is
perpendicular to the normal vector of the face $\mathbf{n}_{f}$ or to the fact
that the integrated Legendre polynomials are evaluated at $1$ or $-1$, which
is zero.
\end{itemize}

Table 2 shows the decomposition of the space $Q_{p+1,p,p}\times Q_{p,p+1,p}%
\times Q_{p,p,p+1}$ for the $\mathbf{\mathcal{H}}(\mathbf{div})$-conforming
hexahedral element.\begin{table}[th]
\begin{center}%
\begin{tabular}
[c]{|c|c|}\hline
Decomposition & Dimension\\\hline
Face-based face functions (lowest order RT) & $6$\\\hline
Face-based face functions (higher order) & $6p(p+2)$\\\hline
interior functions & $3p(p+1)^{2}$\\\hline
Total & $3(p+2)(p+1)^{2}=\dim Q_{p+1,p,p}\times Q_{p,p+1,p}\times Q_{p,p,p+1}%
$\\\hline
\end{tabular}
\end{center}
\caption{Decomposition of the space $\left(  \mathbb{Q}_{n}(K)\right)  ^{3}$
for the $\mathbf{\mathcal{H}}(\mathbf{div})$-conforming hexahedral element.}%
\end{table}

\section{Basis functions for the triangular element}

\label{sec:constructtri2d} The result on the basis construction has been
reported in \cite{Xinpre12}. For the completeness of the current study, we
record the basis functions in this section.

Any point in the 2-simplex $K^{2}$ is uniquely located in terms of the local
coordinate system $(\xi,\eta)$. The vertexes are numbered as $\mathbf{v}%
_{0}(0,0),\mathbf{v}_{1}(1,0),\mathbf{v}_{2}(0,1)$. The barycentric
coordinates are given as
\begin{equation}
\lambda_{0}:=1-\xi-\eta,\quad\lambda_{1}:=\xi,\quad\lambda_{2}:=\eta.
\label{eqn:barycentric2}%
\end{equation}
The directed tangent on a generic edge $\mathbf{e}_{j}=[j_{1},j_{2}]$ is
similarly defined as in (\ref{eqn:edgedef}) for the three-dimensional case. In
the same manner the edge is also parameterized as in (\ref{eqn:edgepara}). A
generic edge can be uniquely identified with
\begin{equation}
\mathbf{e}_{j}:=[j_{1},j_{2}],\quad j_{1}=\{0,1\},\quad j_{1}<j_{2}\leq2,\quad
j=j_{1}+j_{2}. \label{eqn:edgeidentify2}%
\end{equation}
The two-dimensional vectorial curl operator of a scalar quantity, which is
used in our construction, needs a proper definition. We use the one given in
the book \cite{Ray}, \emph{viz}.
\begin{equation}
\mathbf{curl}(u):=\nabla\times u:=\left[  \frac{\partial u}{\partial\eta
},\,-\frac{\partial u}{\partial\xi}\right]  ^{\tau} \label{eqn:scalarveccurl}%
\end{equation}
Based upon the newly created shape functions for the three-dimensional
$\mathbf{\mathcal{H}}(\mathbf{div})$-conforming tetrahedral elements and using
the technique of \emph{dimension reduction} we construct the basis for the
$\mathbf{\mathcal{H}}(\mathbf{div})$-conforming triangular elements in two
dimensions. However, it is easy to see that the two groups for the face
functions cannot be appropriately modified for our purpose. Instead we borrow
the idea of Zaglmayr in the dissertation \cite{Zag}, \emph{viz}., we combine
the edge-based shape functions in \cite{Zag} with our newly constructed
edge-based and bubble interior functions. In \cite{Zag} Zaglmayr had applied
the so-called scaled integrated Legendre polynomials in the construction,
\emph{viz}.
\begin{equation}
\mathcal{L}_{n}^{s}(x,t):=t^{n-1}\int\limits_{-t}^{x}\ell_{n-1}\left(
\frac{\xi}{t}\right)  d\xi,\quad n\geq2,\quad t\in(0,1],
\label{eqn:scalarveccurl1}%
\end{equation}
where $\ell_{n}\left(  x\right)  $ is the n-th order Legender polynomial.

\subsection{Edge functions}

\label{sec:edge2D} For the completeness of our basis construction, in this
subsection we record the results in \cite{Zag}. Associated with each edge the
formulas for these functions are given as
%

\begin{equation}
\Phi_{\mathbf{e}[k_{1},k_{2}]}^{N_{0}} = \lambda_{k_{2}} \nabla\times
\lambda_{k_{1}} - \lambda_{k_{1}} \nabla\times\lambda_{k_{2}}%
\end{equation}
for the lowest-order approximation and
\begin{equation}
\label{eqn:edgehigher}\Phi_{\mathbf{e}[k_{1},k_{2}]}^{j} = \nabla\times\left(
\mathcal{L}_{j+2}^{s} \left(  \gamma_{\mathbf{e}_{k}}, \lambda_{k_{2}} +
\lambda_{k_{1}}\right)  \right)  , \quad j = 0, \cdots, p-1
\end{equation}
for higher-order approximations.

\subsection{Interior functions}

The interior functions are further classified into two categories: edge-based
and bubble interior functions. By construction the normal component of each
interior function vanishes on either edge of the reference 2-simplex $K^{2}$,
\emph{viz}.%

\begin{equation}
\label{eqn:interfun112d}\mathbf{n}^{\mathbf{e}_{j}} \cdot\Phi^{\mathbf{t}} =
0, \quad j = \{1,2,3\},
\end{equation}
where $\mathbf{n}^{\mathbf{e}_{j}}$ is the unit outward normal vector to edge
$\mathbf{e}_{j}$.

\bigskip

\noindent{\underline{Edge-based interior functions}:}

\bigskip

The tangential component of each edge-based function does not vanish on the
associated only edge $\mathbf{e}_{k}:=[k_{1},k_{2}]$ but vanishes the other
two edges, \emph{viz}.
\begin{equation}
\tau^{\mathbf{e}_{j}}\cdot\Phi_{\mathbf{e}[k_{1},k_{2}]}^{\mathbf{t}%
,i}=0,\quad\mathbf{e}_{j}\neq\mathbf{e}_{k}, \label{eqn:edgeinterr112d}%
\end{equation}
where $\tau^{\mathbf{e}_{j}}$ is the directed tangent along the edge
$\mathbf{e}_{j}:=[j_{1},j_{2}]$. The following basis functions are proposed
here as

%

\begin{equation}
\Phi_{\mathbf{e}[k_{1},k_{2}]}^{\mathbf{t},i} = C_{i} \lambda_{k_{1}}
\lambda_{k_{2}} (1 - \lambda_{k_{1}})^{i} P_{i}^{(0,2)}\left(  \frac
{2\lambda_{k_{2}}}{1 - \lambda_{k_{1}}} - 1\right)  \frac{\tau^{\mathbf{e}%
_{k}}}{|\tau^{\mathbf{e}_{k}}|},
\end{equation}
where the function $P_{i}^{(0,2)}(\bullet)$ is the classical
\emph{un-normalized} Jacobi polynomial of degree $i$ with a single variable
\cite{MOS66}, and the scaling coefficient is given by
\begin{equation}
\label{eqn:edge1face22d}C_{i} = \sqrt{2(i + 2)(i + 3)(2i + 3)(2i + 5)}, \quad
i = 0, 1, \cdots, p-2.
\end{equation}

The following orthonormal property of edge-based interior functions can be proved%

\begin{equation}
\label{eqn:edge1faceortho2d}<\Phi_{\mathbf{e}[k_{1},k_{2}]}^{\mathbf{t},m},
\Phi_{\mathbf{e}[k_{1},k_{2}]}^{\mathbf{t},n}>|_{K^{2}} = \delta_{mn},
\quad\{m,n\} = 0,1,\cdots,p-2.
\end{equation}

\bigskip

\noindent{\underline{Interior bubble functions}:}

\bigskip

The interior bubble functions vanish on the entire boundary $\partial K^{2}$
of the reference 2-simplex $K^{2}$. The formulas of these functions are given
as
%

\begin{equation}
\Phi_{m,n}^{\mathbf{t},\vec{e}_{i}}=C_{m,n}\lambda_{0}\lambda_{1}\lambda
_{2}(1-\lambda_{0})^{m}P_{m}^{(2,2)}\left(  \frac{\lambda_{1}-\lambda_{2}%
}{1-\lambda_{0}}\right)  P_{n}^{(2m+5,2)}\left(  2\lambda_{0}-1\right)
\vec{e}_{i},\quad i=1,2, \label{eqn:interdbubtri}%
\end{equation}
where
\[
C_{m,n}=\sqrt{\frac{(m+3)(m+4)(2m+5)(2m+n+6)(2m+n+7)(2m+2n+8)}%
{(m+1)(m+2)(n+1)(n+2)}},
\]
and
\[
0\leq\{m,n\},m+n\leq p-3.
\]

One can again prove the orthonormal property of the interior bubble functions

%

\begin{equation}
<\Phi_{m_{1},n_{1}}^{\mathbf{t},\vec{e}_{i}},\Phi_{m_{2},n_{2}}^{\mathbf{t}%
,\vec{e}_{j}}>|_{K^{2}}=\delta_{m_{1}m_{2}}\delta_{n_{1}n_{2}},
\end{equation}
where
\[
0\leq\{m_{1},m_{2},n_{1},n_{2}\},m_{1}+n_{1},m_{2}+n_{2}\leq
p-3,\,\{i,j\}=1,2.
\]

Table 3 shows the decomposition of the space $\left(  \mathbb{P}%
_{p}(K)\right)  ^{2}$ for the $\mathbf{\mathcal{H}}(\mathbf{div})$-conforming
triangular element.

\begin{table}[th]
\begin{center}%
\begin{tabular}
[c]{|c|c|}\hline
Decomposition & Dimension\\\hline
Edge functions & $3(p+1)$\\\hline
Edge-based interior functions & $3(p-1)$\\\hline
Interior bubble functions & $(p-2)(p-1)$\\\hline
Total & $(p+1)(p+2)=\dim\left(  P_{p}(K)\right)  ^{2}$\\\hline
\end{tabular}
\end{center}
\caption{Decomposition of the space $\left(  \mathbb{P}_{p}(K)\right)  ^{2}$
for the $\mathbf{\mathcal{H}}(\mathbf{div})$-conforming triangular element.}%
\label{tab:deptri}%
\end{table}

\section{Basis functions for the tetrahedral element}

\label{sec:constructtet3} Our constructions are motivated by the work on the
construction of $\mathbf{\mathcal{H}}(\mathbf{div})$-conforming hierarchical
bases for tetrahedral elements \cite{AC03B}. We construct shape functions for
the $\mathbf{\mathcal{H}}(\mathbf{div})$-conforming tetrahedral element on the
canonical reference 3-simplex. The shape functions are grouped into several
categories based upon their geometrical entities on the reference 3-simplex
\cite{AC03B}. The basis functions in each category are constructed so that
they are orthonormal on the reference element.

Any point in the 3-simplex $K^{3}$ is uniquely located in terms of the local
coordinate system $(\xi,\eta,\zeta)$. The vertexes are numbered as
$\mathbf{v}_{0}(0,0,0),\mathbf{v}_{1}(1,0,0),\mathbf{v}_{2}(0,1,0),\mathbf{v}%
_{3}(0,0,1)$. The barycentric coordinates are given as
\begin{equation}
\lambda_{0}:=1-\xi-\eta-\zeta,\quad\lambda_{1}:=\xi,\quad\lambda_{2}%
:=\eta,\quad\lambda_{3}:=\zeta. \label{eqn:barycentric}%
\end{equation}
The directed tangent on a generic edge $\mathbf{e}_{j}=[j_{1},j_{2}]$ is
defined as
\begin{equation}
\tau^{\mathbf{e}_{j}}:=\tau^{\lbrack j_{1},j_{2}]}=\mathbf{v}_{j_{2}%
}-\mathbf{v}_{j_{1}},\quad j_{1}<j_{2}. \label{eqn:edgedef}%
\end{equation}
The edge is parameterized as
\begin{equation}
\gamma_{\mathbf{e}_{j}}:=\lambda_{j_{2}}-\lambda_{j_{1}},\quad j_{1}<j_{2}.
\label{eqn:edgepara}%
\end{equation}
A generic edge can be uniquely identified with
\begin{equation}
\mathbf{e}_{j}:=[j_{1},j_{2}],\quad j_{1}=0,1,2,\quad j_{1}<j_{2}\leq3,\quad
j=j_{1}+j_{2}+\mathrm{sign}(j_{1}), \label{eqn:edgeidentify}%
\end{equation}
where $\mathrm{sign}(0)=0$. Each face on the 3-simplex can be identified by
the associated three vertexes, and is uniquely defined as
\begin{equation}
\mathbf{f}_{j_{1}}:=[j_{2},j_{3},j_{4}],\quad0\leq\{j_{1},j_{2},j_{3}%
,j_{4}\}\leq3,\quad j_{2}<j_{3}<j_{4}. \label{eqn:facedefine}%
\end{equation}

\bigskip

The standard bases in $\mathbb{R}^{n}$ are noted as $\vec{e}_{i}$,
$i=1,\cdots,n$, and $n=\{2,3\}$.

\subsection{Face functions}

\label{sec:face1}

The face functions are further grouped into two categories: edge-based face
functions and face bubble functions.

\bigskip

\noindent{\underline{Edge-based face functions}: $\ $}

\bigskip These functions are associated with the three edges of a certain face
$\mathbf{f}_{j_{1}}$, and by construction all have non-zero normal components
only on the associated face $\mathbf{f}_{j_{1}}$, \emph{viz}.
\begin{equation}
\mathbf{n}^{\mathbf{f}_{j_{k}}}\cdot\Phi_{\mathbf{e}[k_{1},k_{2}]}%
^{\mathbf{f}_{j_{1}},i}=0,\quad j_{k}\neq j_{1}, \label{eqn:edgeface11a}%
\end{equation}
where $\mathbf{n}^{\mathbf{f}_{j_{k}}}$ is the unit outward normal vector to
face $\mathbf{f}_{j_{k}}$.

The edge-based face functions for higher order have been proposed in
\cite{AC03B} as follows:%

\begin{equation}
\widetilde{\Phi}_{\mathbf{e}[k_{1},k_{2}]}^{\mathbf{f}_{j_{1}},i}=l_{i}%
(\gamma_{\mathbf{e}_{k}})\lambda_{k_{1}}\nabla\lambda_{k_{2}}\times
\nabla\lambda_{k_{3}},\;i=0,\cdots,p-1. \label{ac1}%
\end{equation}
For instance, for the face opposite to the vertex $\mathbf{v}_{0}(0,0,0),$
$\mathbf{f}_{0}:=[1,2,3],$ the face functions related to edge $\mathbf{e}%
[1,2]$ are given by%

\begin{equation}
\widetilde{\Phi}_{\mathbf{e}[1,2]}^{\mathbf{f}_{0},i}=l_{i}(\lambda
_{3}-\lambda_{2})\lambda_{1}\nabla\lambda_{2}\times\nabla\lambda
_{3},\;i=0,\cdots,p-1.
\end{equation}

However, it can be checked that the basis function given in (\ref{ac1})
in fact are not independent for $p=2$ (as the sum of the 12 basis functions
given in (\ref{ac1}) for $p=2$ on all 6 faces in fact equals to zero, which
can be easily verified by the symbolic Maple program and the Maple code can be
available from the first author) and thus the proposed basis function is not
complete. \ To remedy this degeneracy, two kinds of constructions of
hierarchical high-order independent edge-based face functions will be
presented here, for which the first one was first reported in \cite{Xinpre12}
while the second kind is proposed here below in (\ref{eqn:edgeface1var}).

\begin{itemize}
\item First kind high-order independent edge-based face functions
\end{itemize}

In \cite{Xinpre12}, the following orthonormal basis functions are given as
%

\begin{equation}
\Phi_{\mathbf{e}[k_{1},k_{2}]}^{\mathbf{f}_{j_{1}},i}=C_{i}\lambda_{k_{3}%
}(1-\lambda_{k_{1}})^{i}P_{i}^{(3,0)}\left(  \frac{2\lambda_{k_{2}}}%
{1-\lambda_{k_{1}}}-1\right)  \frac{\nabla\lambda_{k_{1}}\times\nabla
\lambda_{k_{2}}}{|\nabla\lambda_{k_{1}}\times\nabla\lambda_{k_{2}}|},
\end{equation}
where
\[
C_{i}=\sqrt{3(2i+4)(2i+5)},\quad i=0,1,\cdots,p-1,
\]
and
\[
k_{1}=\{j_{2},j_{3}\},\quad k_{2}=\{j_{3},j_{4}\},\quad k_{1}<k_{2},\quad
k_{3}=\{j_{2},j_{3},j_{4}\}\setminus\{k_{1},k_{2}\}.
\]

One can prove the orthonormal property of these edge-based face functions%

\begin{equation}
<\Phi_{\mathbf{e}[k_{1},k_{2}]}^{\mathbf{f}_{j_{1}},m},\Phi_{\mathbf{e}%
[k_{1},k_{2}]}^{\mathbf{f}_{j_{1}},n}>|_{K^{3}}=\delta_{mn},\quad
\{m,n\}=0,1,\cdots,p-1, \label{eqn:edgefaceortho}%
\end{equation}
where $\delta_{mn}$ is the Kronecker delta. Note that with this construction,
the edge-based face functions are all linearly independent, which is also
verified by the fact that in the spectrum of the mass (Gram) matrix, none of
the eigenvalues is zero.

\bigskip

\begin{itemize}
\item \noindent Second kind high-order independent edge-based face functions
\end{itemize}

\bigskip

An alternative approach using the idea of recursion from \cite{AC03B} can also
be used to construct independent edge-based face functions as follows.

For $p=1$, for each edge we have one face function for this edge as proposed
in \cite{AC03B}%

\begin{equation}
\widetilde{\Phi}_{\mathbf{e}[k_{1},k_{2}]}^{\mathbf{f}_{j_{1}},0}%
=\lambda_{k_{1}}\nabla\lambda_{k_{2}}\times\nabla\lambda_{k_{3}},
\end{equation}
and for $p=2$, one additional new basis function can be constructed as%

\begin{equation}
\widetilde{\Phi}_{\mathbf{e}[k_{1},k_{2}]}^{\mathbf{f}_{j_{1}},1}%
=\lambda_{k_{1}}\lambda_{k_{2}}\nabla\lambda_{k_{3}}\times\nabla\lambda
_{k_{1}},
\end{equation}
which can be shown to satisfy the condition (\ref{eqn:edgeface11a}), and and
for $p\geq3,$ the basis functions are given by
\begin{align}
\widetilde{\Phi}_{\mathbf{e}[k_{1},k_{2}]}^{\mathbf{f}_{j_{1}},i+1}  &
\equiv\ell_{i}(\gamma_{\mathbf{e}_{k}})\widetilde{\Phi}_{\mathbf{e}%
[k_{1},k_{2}]}^{\mathbf{f}_{j_{1}},1}+\ell_{i-1}(\gamma_{\mathbf{e}_{k}%
})\widetilde{\Phi}_{\mathbf{e}[k_{1},k_{2}]}^{\mathbf{f}_{j_{1}}%
,0}\label{eqn:edgeface1var}\\
&  =\ell_{i}(\gamma_{\mathbf{e}_{k}})\left[  \lambda_{k_{1}}\lambda_{k_{2}%
}\nabla\lambda_{k_{3}}\times\nabla\lambda_{k_{1}}\right]  +\ell_{i-1}%
(\gamma_{\mathbf{e}_{k}})\left[  \lambda_{k_{1}}\nabla\lambda_{k_{2}}%
\times\nabla\lambda_{k_{3}}\right]  ,\;i=1,\cdots,p-2.\nonumber
\end{align}
It can be shown again numerically that there are exactly $p$ functions that
are independent and only whose normal component is non-zero only on the
associated edge $\mathbf{e}_{k}$.

\bigskip

\noindent{\underline{Face bubble functions}:}

\bigskip

The face bubble functions which belong to each specific group are associated
with a particular face $\mathbf{f}_{j_{1}}$. They vanish on all edges of the
reference 3-simplex $K^{3}$, and the normal components of which vanish on
other three faces, \emph{viz}.
\begin{equation}
\label{eqn:bubbleface11}\mathbf{n}^{\mathbf{f}_{j_{k}}} \cdot\Phi
_{m,n}^{\mathbf{f}_{j_{1}}} = 0, \quad j_{k} \neq j_{1}.
\end{equation}
The explicit formula is given as
\begin{equation}
\label{eqn:facebub1A}\Phi_{m,n}^{\mathbf{f}_{j_{1}}} = \iota(1 -
\lambda_{j_{2}})^{m} (1 - \lambda_{j_{2}} - \lambda_{j_{3}})^{n}
P_{m}^{(2n+3,2)}\left(  \frac{2\lambda_{j_{3}}}{1 - \lambda_{j_{2}}} -
1\right)  P_{n}^{(0,2)}\left(  \frac{2\lambda_{j_{4}}}{1 - \lambda_{j_{2}} -
\lambda_{j_{3}}} - 1\right)  \frac{\nabla\lambda_{j_{3}} \times\nabla
\lambda_{j_{4}}}{|\nabla\lambda_{j_{3}} \times\nabla\lambda_{j_{4}}|},
\end{equation}
where
\begin{equation}
\label{eqn:facebub2}\iota= C_{m}^{n} \lambda_{j_{2}} \lambda_{j_{3}}
\lambda_{j_{4}},
\end{equation}
where
\begin{equation}
\label{eqn:facebub3}C_{m}^{n} = \frac{\sqrt{(2n + 3)(m + n + 3)(m + 2n + 4)(m
+ 2n + 5)(2m + 2n + 7)(2m + 2n + 8)(2m + 2n + 9)}} {\sqrt{(m + 1)(m + 2)}},
\end{equation}
and
\begin{equation}
\label{eqn:facebub4}0 \le\{m,n\}, m+n \le p-3.
\end{equation}
By construction the face bubble functions share again the orthonormal property
on the reference 3-simplex $K^{3}$:
\begin{equation}
\label{eqn:facebubbleorthonormal1a}<\Phi_{m_{1},n_{1}}^{\mathbf{f}_{j_{1}}},
\Phi_{m_{2},n_{2}}^{\mathbf{f}_{j_{1}}}>|_{K^{3}} = \delta_{m_{1} m_{2}}%
\delta_{n_{1} n_{2}}, \quad0 \le\{m_{1}, m_{2}, n_{1}, n_{2}\}, m_{1} + n_{1},
m_{2} + n_{2} \le p-3.
\end{equation}

\subsection{Interior functions}

\label{sec:interior1} The interior functions are classified into three
categories: edge-based, face-based and bubble interior functions. By
construction the normal component of each interior function vanishes on either
face of the reference 3-simplex $K^{3}$, \emph{viz}.
\begin{equation}
\label{eqn:interfun11}\mathbf{n}^{\mathbf{f}_{j}} \cdot\Phi^{\mathbf{t}} = 0,
\quad j = \{0,1,2,3\}.
\end{equation}

\bigskip

\noindent{\underline{Edge-based interior functions}:}

\bigskip

The tangential component of each edge-based function does not vanish on the
associated only edge $\mathbf{e}_{k} := [k_{1}, k_{2}]$ but vanishes all other
five edges, \emph{viz}.
\begin{equation}
\label{eqn:edgeinterr11}\tau^{\mathbf{e}_{j}} \cdot\Phi_{\mathbf{e}%
[k_{1},k_{2}]}^{\mathbf{t},i} = 0, \quad\mathbf{e}_{j} \neq\mathbf{e}_{k},
\end{equation}
where $\tau^{\mathbf{e}_{j}}$ is the directed tangent along the edge
$\mathbf{e}_{j} := [j_{1},j_{2}]$. The shape functions are given as

%

\begin{equation}
\Phi_{\mathbf{e}[k_{1},k_{2}]}^{\mathbf{t},i}=C_{i}\lambda_{k_{1}}%
\lambda_{k_{2}}(1-\lambda_{k_{1}})^{i}P_{i}^{(1,2)}\left(  \frac
{2\lambda_{k_{2}}}{1-\lambda_{k_{1}}}-1\right)  \frac{\tau^{\mathbf{e}_{k}}%
}{|\tau^{\mathbf{e}_{k}}|},
\end{equation}
where
\[
C_{i}=(i+3)\sqrt{\frac{(2i+4)(2i+5)(2i+7)}{i+1}},\quad i=0,1,\cdots,p-2.
\]

Again one can prove the orthonormal property of edge-based interior functions:%

\begin{equation}
\label{eqn:edge1faceortho}<\Phi_{\mathbf{e}[k_{1},k_{2}]}^{\mathbf{t},m},
\Phi_{\mathbf{e}[k_{1},k_{2}]}^{\mathbf{t},n}>|_{K^{3}} = \delta_{mn},
\quad\{m,n\} = 0,1,\cdots,p-2.
\end{equation}

\bigskip

\noindent{\underline{Face-based interior functions}:}

\bigskip

These functions which are associated with a particular face $\mathbf{f}%
_{j_{1}}$ have non-zero tangential components on their associated face only,
and have no contribution to the tangential components on all other three
faces, \emph{viz}.
\begin{equation}
\mathbf{n}^{\mathbf{f}_{j_{k}}}\times\Phi_{m,n}^{\mathbf{t},\mathbf{f}_{j_{1}%
}}=\mathbf{0},\quad j_{k}\neq j_{1}. \label{eqn:inter1face11}%
\end{equation}
Further each face-based interior function vanishes on all the edges of the
3-simplex $K^{3}$, \emph{viz}.
\begin{equation}
\tau^{\mathbf{e}_{k}}\cdot\Phi_{m,n}^{\mathbf{t},\mathbf{f}_{j_{1}}}=0.
\label{eqn:face22interr11}%
\end{equation}
The formulas of these functions are given as
%

\[
\Phi_{m,n}^{\mathbf{t},\mathbf{f}_{j_{1}}^{1}}=\iota(1-\lambda_{j_{2}}%
)^{m}(1-\lambda_{j_{2}}-\lambda_{j_{3}})^{n}P_{m}^{(2n+3,2)}\left(
\frac{2\lambda_{j_{3}}}{1-\lambda_{j_{2}}}-1\right)  P_{n}^{(0,2)}\left(
\frac{2\lambda_{j_{4}}}{1-\lambda_{j_{2}}-\lambda_{j_{3}}}-1\right)
\frac{\tau^{\lbrack j_{2},j_{3}]}}{\left\vert \tau^{\lbrack j_{2},j_{3}%
]}\right\vert },
\]%
\begin{equation}
\Phi_{m,n}^{\mathbf{t},\mathbf{f}_{j_{1}}^{2}}=\iota(1-\lambda_{j_{2}}%
)^{m}(1-\lambda_{j_{2}}-\lambda_{j_{3}})^{n}P_{m}^{(2n+3,2)}\left(
\frac{2\lambda_{j_{3}}}{1-\lambda_{j_{2}}}-1\right)  P_{n}^{(0,2)}\left(
\frac{2\lambda_{j_{4}}}{1-\lambda_{j_{2}}-\lambda_{j_{3}}}-1\right)
\frac{\tau^{\lbrack j_{2},j_{4}]}}{\left\vert \tau^{\lbrack j_{2},j_{4}%
]}\right\vert }, \label{eqn:faceinter22}%
\end{equation}
where $\iota$ is given in (\ref{eqn:facebub2}) and $0\leq\{m,n\},m+n\leq p-3$.
The face-based interior functions enjoy the orthonormal property on the
reference 3-simplex $K^{3}$:%

\begin{equation}
\label{eqn:faceinterorthonormal}<\Phi_{m_{1},n_{1}}^{\mathbf{t},
\mathbf{f}_{j_{1}}^{i}},\Phi_{m_{2},n_{2}}^{\mathbf{t}, \mathbf{f}_{j_{1}}%
^{i}}>|_{K^{3}} = \delta_{m_{1} m_{2}}\delta_{n_{1} n_{2}}, i = \{1,2\}, 0
\le\{m_{1}, m_{2}, n_{1}, n_{2}\}, m_{1} + n_{1}, m_{2} + n_{2} \le p-3.
\end{equation}

\bigskip

\noindent{\underline{Interior bubble functions}:}

\bigskip

The interior bubble functions vanish on the entire boundary $\partial K^{3}$
of the reference 3-simplex $K^{3}$. The formulas of these functions are given
as
%

\begin{equation}
\Phi_{\ell,m,n}^{\mathbf{t},\vec{e}_{i}}=\chi P_{\ell}^{(2m+2n+8,2)}\left(
2\lambda_{1}-1\right)  P_{m}^{(2n+5,2)}\left(  \frac{2\lambda_{2}}%
{1-\lambda_{1}}-1\right)  P_{n}^{(2,2)}\left(  \frac{2\lambda_{3}}%
{1-\lambda_{1}-\lambda_{2}}-1\right)  \vec{e}_{i},\,i=1,2,3,
\label{eqn:interbub}%
\end{equation}
where
\[
\chi=C_{\ell,m,n}\lambda_{0}\lambda_{1}\lambda_{2}\lambda_{3}(1-\lambda
_{1})^{m}(1-\lambda_{1}-\lambda_{2})^{n},
\]
where
\[
C_{\ell,m,n}=C_{\ell,m,n}^{1}C_{\ell,m,n}^{2},
\]
where
\[
C_{\ell,m,n}^{1}=\sqrt{\frac{(\ell+2m+2n+9)(\ell+2m+2n+10)(2\ell
+2m+2n+11)(m+2n+6)}{(\ell+1)(m+1)(n+1)}},
\]%
\[
C_{\ell,m,n}^{2}=\sqrt{\frac{(m+2n+7)(2m+2n+8)(n+3)(n+4)(2n+5)}{(\ell
+2)(m+2)(n+2)}},
\]
and
\[
0\leq\{\ell,m,n\},\ell+m+n\leq p-4.
\]
Again, one can show the orthonormal property of the interior bubble functions

%

\[
<\Phi_{\ell_{1},m_{1},n_{1}}^{\mathbf{t},\vec{e}_{i}},\Phi_{\ell_{2}%
,m_{2},n_{2}}^{\mathbf{t},\vec{e}_{j}}>|_{K^{3}}=\delta_{\ell_{1}\ell_{2}%
}\delta_{m_{1}m_{2}}\delta_{n_{1}n_{2}},
\]
where
\[
0\leq\{\ell_{1},\ell_{2},m_{1},m_{2},n_{1},n_{2}\},\ell_{1}+m_{1}+n_{1}%
,\ell_{2}+m_{2}+n_{2}\leq p-4,\,\{i,j\}=1,2,3.
\]

In Table 4 we summarize the decomposition of the space $\left(  \mathbb{P}%
_{p}(K)\right)  ^{3}$ for the $\mathbf{\mathcal{H}}(\mathbf{div})$-conforming
tetrahedral element.

\begin{table}[th]
\begin{center}%
\begin{tabular}
[c]{|c|c|}\hline
Decomposition & Dimension\\\hline
Edge-based face functions & $12p$\\\hline
Face bubble functions & $2(p-2)(p-1)$\\\hline
Edge-based interior functions & $6(p-1)$\\\hline
Face-based interior functions & $4(p-2)(p-1)$\\\hline
Interior bubble functions & $(p-3)(p-2)(p-1)/2$\\\hline
Total & $(p+1)(p+2)(p+3)/2=\dim\left(  P_{p}(K)\right)  ^{3}$\\\hline
\end{tabular}
\end{center}
\caption{Decomposition of the space $\left(  \mathbb{P}_{p}(K)\right)  ^{3}$
for the $\mathbf{\mathcal{H}}(\mathbf{div})$-conforming tetrahedral element.}%
\label{tab:deptet}%
\end{table}

\section{The divergence-free condition}

\label{sec:divcond} For a $p$-th order polynomial approximation and to ensure
the divergence-free condition, the idea is to include a higher-order interior
bubble function $\chi_{b}$, say $\chi_{b}^{p+1}$ or $\chi_{b}^{q}%
,q=\max(m,p+1)$ for the simplexes ( $m=3$ for triangular element and $m=4$ for
tetrahedral element) and $q=\max(2,p+1)$ for rectangular and hexahedral
elements , as an extra basis function. To satisfy the condition
\begin{equation}
\nabla\cdot\mathbf{B}=0 \label{eqn:conddivB}%
\end{equation}
for the Maxwell equation or
\begin{equation}
\nabla\cdot\mathbf{u}=0 \label{eqn:conddivvelocity}%
\end{equation}
for the incompressible fluid flow, one can impose the condition
\begin{equation}
\nabla\cdot(C_{\chi}\chi_{b}^{p+1})=0 \label{eqn:condimpose}%
\end{equation}
to solve the unknown coefficient $C_{\chi}$. For the triangular element, the
interior bubble function for degree $p$ is given in equation
(\ref{eqn:interdbubtri}). For the tetrahedral element, the interior bubble
function of degree $p$ is given in equation (\ref{eqn:interbub}). For the
quadrilateral element, one can construct the interior bubble function of
degree $p$ as
\begin{equation}
\chi_{\mathcal{Q}_{b}}^{p}=L_{p+2}(2\xi-1)L_{2}(2\eta-1)\hat{\mathbf{i}}_{\xi
}+L_{2}(2\xi-1)L_{p+2}(2\eta-1)\hat{\mathbf{j}}_{\eta}.
\label{eqn:quadinterbubble}%
\end{equation}
For the hexahedral element, one can construct the interior bubble function of
degree $p$ as
\begin{equation}
\chi_{\mathcal{H}_{b}}^{p}=L_{p+2}(\xi_{1})L_{2}(\eta_{2})L_{2}(\zeta_{3}%
)\hat{\mathbf{i}}_{\xi}+L_{2}(\xi_{1})L_{p+2}(\eta_{2})L_{2}(\zeta_{3}%
)\hat{\mathbf{j}}_{\eta}+L_{2}(\xi_{1})L_{2}(\eta_{2})L_{p+2}(\zeta_{3}%
)\hat{\mathbf{k}}_{\zeta}. \label{eqn:hexainterbubble}%
\end{equation}

\section{Conditioning of matrices}

The purpose of this section is twofold. Firstly, we check numerically that the
newly constructed basis functions for $\mathbf{\mathcal{H}}(\mathbf{div}%
)$-conforming triangular and tetrahedral elements are linearly independent,
which is manifested by the fact that for each particular approximation order
up to degree four, the condition number of the corresponding mass matrix is
finite. Secondly, we want to show that for the approximation up to order three
both the mass and stiffness matrices are reasonably well-conditioned.

The components of the mass matrix are defined as
\begin{equation}
\label{eqn:components1}M_{\ell_{1},\ell_{2}}:=\,<\Phi_{\ell_{1}},\Phi
_{\ell_{2}}>|_{K^{d}}, \quad d=2,3.
\end{equation}
The mass matrix $M$ is symmetric and positive definite, and therefore has real
positive eigenvalues. The condition number of a real symmetric positive
definite matrix $A$ is calculated by the formula
\begin{equation}
\label{eqn:condition}\kappa(A) = \frac{\lambda_{\max}}{\lambda_{\min}},
\end{equation}
where $\lambda_{\max}$ and $\lambda_{\min}$ are the maximum and minimum
eigenvalues of the matrix $A$, respectively. For the incompressible fluid
flows, e.g., governed by the Navier-Stokes equations \cite{Kan} or by the
magnetohydrodynamics equations \cite{Dan} the authors \cite{Kan, Dan} have
applied the mixed finite element for the spatial discretization. In
particular, they \cite{Kan, Dan} have used the $\mathbf{\mathcal{H}%
}(\mathbf{div})$-conforming element for the Laplacian $\Delta\mathbf{u}$ of
the velocity $\mathbf{u}$. In this case, we have the stiffness matrix $S$,
which is defined component-wise as
\begin{equation}
\label{eqn:components2}S_{\ell_{1},\ell_{2}}:=\,<\nabla\Phi_{\ell_{1}} :
\nabla\Phi_{\ell_{2}}>|_{K^{d}}, \quad d=2,3.
\end{equation}
The stiffness matrix $S$ is symmetric and semi-positive definite, and
therefore has real non-negative eigenvalues. The condition number of the
stiffness matrix $S$ is calculated by the formula (\ref{eqn:condition}) with
the zero eigenvalue excluded.

With the triangular element and for the polynomial approximations $p = \{1, 2,
3, 4\}$, the conditioning is summarized in Table \ref{tab:condtri}.
\begin{table}[th]
\begin{center}%
\begin{tabular}
[c]{|c|c|c|}\hline
Order $p$ & Mass & Stiffness\\\hline
1 & 2.016e1 & 1.040e1\\\hline
2 & 8.804e1 & 5.959e1\\\hline
3 & 9.847e2 & 4.197e2\\\hline
4 & 1.286e4 & 8.843e3\\\hline
\end{tabular}
\end{center}
\caption{Condition numbers of the mass matrix $M$ and stiffness matrix $S$
from the basis for the $\mathbf{\mathcal{H}}(\mathbf{div})$-conforming
triangular element.}%
\label{tab:condtri}%
\end{table}From the table we can see that the condition number is bounded for
each order of approximation. Moreover, up to the fourth order, the mass and
stiffness matrices are both well conditioned.

With the tetrahedral element and for the polynomial approximations $p = \{1,
2, 3, 4\}$, the condition numbers of the mass matrix are shown in Table
\ref{tab:condtet}.

\begin{table}[th]
\begin{center}%
\begin{tabular}
[c]{|c|c|c|c|}\hline
\multicolumn{1}{|c|}{Order} & Mass & Stiffness & Ratio\\\cline{2-4}%
$p$ & First kind \vline\ Second kind & First kind \vline\ Second kind & Mass
\vline\ Stiff.\\\hline
1 & 3.084e1 \vline\ 3.084e1 & 1.989e1 \vline\ 1.989e1 & 1.000e0
\vline\ 1.000e0\\\hline
2 & 6.987e3 \vline\ 7.733e4 & 3.395e3 \vline\ 5.917e4 & 0.090e0
\vline\ 0.057e0\\\hline
3 & 3.412e6 \vline\ 2.289e6 & 1.094e6 \vline\ 1.191e6 & 1.491e0
\vline\ 0.919e0\\\hline
4 & 5.972e9 \vline\ 2.717e7 & 2.883e9 \vline\ 2.372e7 & 2.198e2
\vline\ 1.215e2\\\hline
\end{tabular}
\end{center}
\caption{Condition numbers of the mass matrix $M$ and stiffness matrix $S$
from the bases with two different kinds of edge-based face functions for the
$\mathbf{\mathcal{H}}(\mathbf{div})$-conforming tetrahedral element.}%
\label{tab:condtet}%
\end{table}Again from this table we see that the condition number is finite
for each order of approximation. Further up to the third order, both the mass
and stiffness matrices are well conditioned. For order $p=2$ the conditioning
is better with the first kind edge-based face basis while for $p=4$, the
conditioning is better with the second kind edge-based face basis. For the
third order $p=3$, the performance with both kinds of edge-based face bases is
about the same.

\section{Concluding remarks}

\label{sec:conclusion} In this paper we focus our attention on hierarchical
$\mathbf{\mathcal{H}}(\mathbf{div})$ basis functions for solving the
magnetohydrodynamics (MHD) equations numerically so that the divergence-free
condition on the magnetic field is rigorously guaranteed. The idea is to use
an interior bubble function from the proposed high-order hierarchical basis as
the additional freedom to impose the divergence-free Gauge condition for the
magnetic field. We have summarized four bases for the $\mathbf{\mathcal{H}%
}(\mathbf{div})$-conforming elements, viz. the quadrilateral and triangular
elements for 2-D and the hexahedral and tetrahedral elements for 3-D. The
linear independence of the basis functions for the two simplicial elements has
numerically been checked. Good matrix (mass and stiffness) conditioning has
also been shown up to the fourth order for 2-D and up to the third order for
3-D. Further work will include the implementation of the proposed
divergence-free basis to solve the magnetohydrodynamics (MHD) equations in 2-D
and 3-D.

\section*{Acknowledgment}

This research is supported in part by a DOE grant DEFG0205ER25678 and a NSF grant DMS-1005441.

\end{document}